\magnification 1200

  \input amssym
  \input miniltx
  \input pictex

  %
  \font \bbfive = bbm5
  \font \bbseven = bbm7
  \font \bbten = bbm10
  \font \eightbf = cmbx8
  \font \eighti = cmmi8 \skewchar \eighti = '177
  \font \eightit = cmti8
  \font \eightrm = cmr8
  \font \eightsl = cmsl8
  \font \eightsy = cmsy8 \skewchar \eightsy = '60
  \font \eighttt = cmtt8 \hyphenchar \eighttt = -1

  \font \sixi = cmmi6 \skewchar \sixi = '177
  \font \sixrm = cmr6
  \font \sixsy = cmsy6 \skewchar \sixsy = '60
  \font \tensc = cmcsc10

  \scriptfont \bffam = \bbseven
  \scriptscriptfont \bffam = \bbfive
  \textfont \bffam = \bbten

  \newskip \ttglue

  \def \eightpoint {\def \rm {\fam 0 \eightrm }\relax
  \textfont 0= \eightrm
  \scriptfont 0 = \sixrm \scriptscriptfont 0 = \fiverm
  \textfont 1 = \eighti
  \scriptfont 1 = \sixi \scriptscriptfont 1 = \fivei
  \textfont 2 = \eightsy
  \scriptfont 2 = \sixsy \scriptscriptfont 2 = \fivesy
  \textfont 3 = \tenex
  \scriptfont 3 = \tenex \scriptscriptfont 3 = \tenex
  \def \it {\fam \itfam \eightit }\relax
  \textfont \itfam = \eightit
  \def \sl {\fam \slfam \eightsl }\relax
  \textfont \slfam = \eightsl
  \def \bf {\fam \bffam \eightbf }\relax
  \textfont \bffam = \bbseven
  \scriptfont \bffam = \bbfive
  \scriptscriptfont \bffam = \bbfive
  \def \tt {\fam \ttfam \eighttt }\relax
  \textfont \ttfam = \eighttt
  \tt \ttglue = .5em plus.25em minus.15em
  \normalbaselineskip = 9pt
  \def \MF {{\manual opqr}\-{\manual stuq}}\relax
  \let \sc = \sixrm
  \let \big = \eightbig
  \setbox \strutbox = \hbox {\vrule height7pt depth2pt width0pt}\relax
  \normalbaselines \rm }

  \def \withfont #1#2{\font \auxfont =#1 {\auxfont #2}}

  %

  \def \TRUE {Y}
  \def \FALSE {N}
  \def \EMPTY {}

  \def \ifundef #1{\expandafter \ifx \csname #1\endcsname \relax }

  \def \undefrule{\kern 2pt \vrule width 5pt height 5pt depth 0pt \kern 2pt}
  \def \UndefLabels{}
  \def \possundef #1{\ifundef {#1}\undefrule {\eighttt #1}\undefrule
    \global \edef \UndefLabels{\UndefLabels#1\par }
  \else \csname #1\endcsname \fi }

  %

  \newcount \secno \secno = 0
  \newcount \stno \stno = 0
  \newcount \eqcntr \eqcntr = 0

  \ifundef {showlabel} \global \def \showlabel {\FALSE} \fi  
  \ifundef {auxfile}   \global \def \auxfile   {\TRUE} \fi

  \def \define #1#2{\global \expandafter \edef \csname #1\endcsname {#2}}
  \long \def \error #1{\medskip \noindent {\bf ******* #1}}
  \def \fatal #1{\error{#1\par Exiting...}\end }

  \def \advseqnumbering {\global \advance \stno by 1 \global \eqcntr =0}

  \def \current {\ifnum \secno = 0 \number \stno \else \number \secno \ifnum \stno = 0 \else .\number \stno \fi \fi}

  \begingroup \catcode `\@=0 \catcode `\\=11 @global @def @textbackslash {\} @endgroup
  \begingroup \catcode `\%=11                \global \def \percent       {
  \def \space { }

  %
  \def \deflabel #1#2{%
    \if\TRUE\showlabel \hbox {\sixrm [[ #1 ]]} \fi
    \ifundef {#1PrimarilyDefined}%
      \define{#1}{#2}%
      \define{#1PrimarilyDefined}{#2}%
      \if\TRUE\auxfile \immediate \write 1 {\textbackslash newlabel {#1}{#2}}\fi
    \else
      \edef \old {\csname #1\endcsname}%
      \edef \new {#2}%
      \if \old \new \else \fatal{Duplicate definition for label ``{\tt #1}'', already defined as ``{\tt \old}''.}\fi
      \fi}

  \def \newlabel #1#2{\define{#1}{#2}}  
  \def \label #1 {\deflabel {#1}{\current }}

  \def \equationmark #1 {\ifundef {InsideBlock}
	  \advseqnumbering
	  \eqno {(\current )}
	  \deflabel {#1}{\current }
	\else
	  \global \advance \eqcntr by 1
	  \edef \subeqmarkaux {\current .\number \eqcntr }
	  \eqno {(\subeqmarkaux )}
	  \deflabel {#1}{\subeqmarkaux }
	\fi }

  \def \split #1.#2.#3.#4;{\global \def \parone {#1}\global \def \partwo {#2}\global \def \parthree {#3}\global \def \parfour {#4}}
  \def \NA {NA}
  \def \ref #1{\split #1.NA.NA.NA;(\possundef {\parone }\ifx \partwo \NA \else .\partwo \fi )}
  \def \redundantref #1#2{\ref {#2}}

  %
  \newcount \bibno \bibno = 0
  \def \newbib #1#2{\define{#1}{#2}} 

  \def \Bibitem #1 #2; #3; #4 \par{\smallbreak
    \global \advance \bibno by 1
    \item {[\possundef{#1}]} #2, {``#3''}, #4.\par
    \ifundef {#1PrimarilyDefined}\else
      \fatal{Duplicate definition for bibliography item ``{\tt #1}'', already defined in ``{\tt [\csname #1\endcsname]}''.}
      \fi
	\ifundef {#1}\else
	  \edef \prevNum{\csname #1\endcsname}
	  \ifnum \bibno=\prevNum \else
		\error{Mismatch bibliography item ``{\tt #1}'', defined earlier (in aux file ?) as ``{\tt \prevNum}'' but should be
	``{\tt \number\bibno}''.  Running again should fix this.}
		\fi
	  \fi
    \define{#1PrimarilyDefined}{#2}%
    \if\TRUE\auxfile \immediate\write 1 {\textbackslash newbib {#1}{\number\bibno}}\fi}

  \def \jrn #1, #2 (#3), #4-#5;{{\sl #1}, {\bf #2} (#3), #4--#5}
  \def \Article #1 #2; #3; #4 \par{\Bibitem #1 #2; #3; \jrn #4; \par}

  \def \references {\begingroup \bigbreak \eightpoint \centerline {\tensc References} \nobreak \medskip \frenchspacing }

  %

  \catcode `\@=11
  \def \c@itrk #1{{\bf \possundef {#1}}} 
  \def \c@ite #1{{\rm [\c@itrk{#1}]}}
  \def \sc@ite [#1]#2{{\rm [\c@itrk{#2}\hskip 0.7pt:\hskip 2pt #1]}}
  \def \du@lcite {\if \pe@k [\expandafter \sc@ite \else \expandafter \c@ite \fi }
  \def \cite {\futurelet\pe@k \du@lcite }
  \catcode `\@=12

  %
  \def \Headlines #1#2{\nopagenumbers
    \headline {\ifnum \pageno = 1 \hfil
    \else \ifodd \pageno \tensc \hfil \lcase {#1} \hfil \folio
    \else \tensc \folio \hfil \lcase {#2} \hfil
    \fi \fi }}

  \def \title #1{\medskip\centerline {\withfont {cmbx12}{\ucase{#1}}}}

  \def \Subjclass #1#2{\footnote {\null }{\eightrm #1 \eightsl Mathematics Subject Classification:  \eightrm #2.}}

  \long \def \Quote #1\endQuote {\begingroup \leftskip 35pt \rightskip 35pt
\parindent 17pt \eightpoint #1\par \endgroup }
  \long \def \Abstract #1\endAbstract {\bigskip \Quote \noindent #1\endQuote }
  \def \Address #1#2{\bigskip {\tensc #1 \par \it E-mail address: \tt #2}}
  \def \Authors #1{\bigskip \centerline {\tensc #1}}
  \def \Note #1{\footnote {}{\eightpoint #1}}
  \def \Date #1 {\Note {\it Date: #1.}}

  \def \part #1#2{\vfill \eject \null \vskip 0.3\vsize
    \withfont{cmbx10 scaled 1440}{\centerline{PART #1} \vskip 1.5cm \centerline{#2}} \vfill\eject }

  %

  \def \fix {\smallskip \noindent $\blacktriangleright $\kern 12pt}
  \def \iskip {\medskip\noindent}
  \def \warn {??\vrule height 7pt width 12pt??}
  \def \newpage {\vfill \eject }

  \def \ucase #1{\edef \auxvar {\uppercase {#1}}\auxvar }
  \def \lcase #1{\edef \auxvar {\lowercase {#1}}\auxvar }

  \def \section #1 \par {\global \advance \secno by 1 \stno = 0
    %
    \goodbreak \bigbreak
    \noindent {\bf \number \secno .\enspace #1.}
    \nobreak \medskip \noindent }

  \def \state #1 #2\par {\begingroup \def \InsideBlock {} \medbreak \noindent \advseqnumbering {\bf \current .\enspace
#1.\enspace \sl #2\par }\medbreak \endgroup }

  \def \definition #1\par {\state Definition \rm #1\par }

  \long \def \Proof #1\endProof {\begingroup \def \InsideBlock {} \medbreak \noindent {\it Proof.\enspace }#1
\ifmmode \eqno \endproofmarker $$ \else \hfill $\endproofmarker $ \looseness = -1 \fi \medbreak \endgroup }

  \def \$#1{#1 $$$$ #1}
  \def \explica #1#2{\mathrel {\buildrel \hbox {\sixrm #2} \over #1}}
  \def \explain #1#2{\explica{#1}{\ref{#2}}}  
  \def \=#1{\explain {=}{#1}}

  \def \pilar #1{\vrule height #1 width 0pt}
  \def \stake #1{\vrule depth  #1 width 0pt}

  \newcount \fnctr \fnctr = 0
  \def \fn #1{\global \advance \fnctr by 1
    \edef \footnumb {$^{\number \fnctr }$}%
    \footnote {\footnumb }{\eightpoint #1\par \vskip -10pt}}

  \def \text #1{\hbox {#1}}
  \def \bool #1{[{\scriptstyle #1}]\,}
  \def \equal #1#2{\bool {#1=#2}}

  %
  \def \iItem {\smallskip }
  \def \Item #1{\smallskip \item {{\rm #1}}}
  \newcount \zitemno \zitemno = 0

  \def \izitem {\global \zitemno = 0}
  \def \zitemplus {\global \advance \zitemno by 1 \relax }
  \def \rzitem {\romannumeral \zitemno }
  \def \rzitemplus {\zitemplus \rzitem } 
  \def \zitem {\Item {{\rm (\rzitemplus )}}}
  \def \Zitem {\izitem \zitem }
  \def \zitemmark #1 {\deflabel {#1}{\rzitem }}

  \newcount \nitemno \nitemno = 0
  \def \initem {\nitemno = 0}
  \def \nitem {\global \advance \nitemno by 1 \Item {{\rm (\number \nitemno )}}}

  \newcount \aitemno \aitemno = -1
  \def \boxlet #1{\hbox to 6.5pt{\hfill #1\hfill }}
  \def \iaitem {\aitemno = -1}
  \def \aitemconv {\ifcase \aitemno a\or b\or c\or d\or e\or f\or g\or
h\or i\or j\or k\or l\or m\or n\or o\or p\or q\or r\or s\or t\or u\or
v\or w\or x\or y\or z\else zzz\fi }
  \def \aitem {\global \advance \aitemno by 1\Item {(\boxlet \aitemconv )}}
  \def \aitemmark #1 {\deflabel {#1}{\aitemconv }}

  \def\Bitem{\Item{$\bullet$}}

  \def \Case #1:{\medskip \noindent {\tensc Case #1:}}

  %
  \def \<{\left \langle \vrule width 0pt depth 0pt height 8pt }
  \def \>{\right \rangle }
  \def \({\big (}
  \def \){\big )}
  \def \ds {\displaystyle }
  \def \and {\hbox {,\quad and \quad }}
  \def \calcat #1{\,{\vrule height8pt depth4pt}_{\,#1}}
  \def \imply {\mathrel {\Rightarrow }}
  \def \IFF {\kern 7pt\Leftrightarrow \kern 7pt}
  \def \IMPLY {\kern 7pt \Rightarrow \kern 7pt}
  \def \for #1{,\quad \forall \,#1}
  \def \endproofmarker {\square } 
  \def \"#1{{\it #1}\/} \def \umlaut #1{{\accent "7F #1}}
  \def \inv {^{-1}}
  \def \*{\otimes }
  \def \caldef #1{\global \expandafter \edef \csname #1\endcsname {{\cal #1}}}
  \def \bfdef #1{\global \expandafter \edef \csname #1\endcsname {{\bf #1}}}
  \bfdef N \bfdef Z \bfdef C \bfdef R

  %

  \if \TRUE \auxfile
    \IfFileExists {\jobname.aux}{\input \jobname.aux}{\null}
    \immediate \openout 1 \jobname.aux
    \fi

  \def \close {\if \EMPTY \UndefLabels \else
      \immediate \write 1{\percent\space *** There were undefined labels ***}
      \message {*** There were undefined labels ***} \iskip
      ****************** \ Undefined Labels: \tt \par \UndefLabels
      \fi
    \if \TRUE \auxfile \closeout 1 \fi
    \par \vfill \supereject \end }

  \def \inpatex #1.tex{\input #1.atex}

  %

  \def \Caixa #1{\setbox 1=\hbox {$#1$\kern 1pt}\global \edef \tamcaixa {\the \wd 1}\box 1}
  \def \caixa #1{\hbox to \tamcaixa {$#1$\hfil }}

  \def \med #1{\mathop {\textstyle #1}\limits }
  \def \medoplus {\med \bigoplus }
  \def \medsum {\med \sum }
  \def \medprod {\med \prod }
  \def \medcup {\med \bigcup }

  \def \nopar {\nopagenumbers \parindent 0pt \parskip 10pt}


  %

  \def \src {d}	 \def\sr#1{\src(#1)}               
  \def \ran {r}	 \def\rn#1{\ran(#1)}               
  \def \vr {x}	                                   
  \def \vro {y}	                                   
  \def \ed {e}                                     
  \def \oed {f}	                                   
  \def \eproj {f}                                  
  \def \s {s}  	                                   
  \def \auto {\sigma}                                   
  \def \ts {\tilde \s }                            
  \def \tp {\tilde p}                              
  \def \tu {\tilde u}                              
  \def \calS {{\cal S}}                            
  \def \E {{\cal E}}                               
  \def \SGE {{\cal S}_{G,E}}                       
  \def \EGE {{\cal E}}                             
  \def \SE {{\cal S}_E}                            
  \def \GpdGE {{\cal G}\tight (\SGE )}             

  \def \q {\breve }
  \def \corona {{\q G}}
  \def \lag {\ell}
  \def \O {{\cal O}}

  \def \g {{\bf g}}  
  \def \cyl #1{Z(#1)}

  \def \germ #1#2#3#4{\big [#1,#2,#3;\,#4\big ]}
  \def \trunc #1#2{#1|_{#2}}

  \def \Lin {{\cal L}}
  \def \acite [#1]{\cite [#1]{actions}}



  \def \S {{\cal S}}
  \def \E {{\cal E}}
  \def \interior #1{{\buildrel \circ \over #1}}

  \def \TF {T\!F}
  \def \J {{\cal J}}
  \def \tight {_{\rm tight}}
  \def \red {_{\rm red}}


  \def \G {{\cal G}}
  \def \Dm #1#2{D^{#1}_{#2}}
  \def \D #1#2{D^{#1}\{#2\}}  
  \def \act {\alpha}
  \def \actfil {\beta}
  \def \X #1#2{{\cal F}^{\,#1}_{#2}} 
  \def \Gth {\G \tight (\S )}
  \def \Eth {\hat \E \tight }
  \def \zero {^{(0)}}


  \def \superfixed {fixed}


  \def \orb #1{\hbox {Orb}(#1)}


  \def \clos #1{\overline {#1}}



  \def \Data {G,E,\varphi}
  \def \OGE {{\cal O}_{G,E}}
  \def \OAB {{\cal O}_{A,B}}


  \def \proj {q} 
  \def \gp {v} 
  \def \rep {V} 
  \def \t {t} 
  \def \modmap {\Psi} 
  \def \repalg {\psi} 
  \def \projMod {Q} 


  \def \mathcal #1{{\cal #1}}
  \def \text #1{\hbox {#1}}
  \def \mathbb #1{{\bf #1}}
  \def \mbox #1{\hbox {#1}}
  \def \frac #1#2{{#1 \over #2}}

  \def \vspace #1{\vskip #1}
  \def \hspace #1{\hskip #1}
  \def \lhd {$ is a normal subgroup of $}

  \newcount \itemdpt \itemdpt = 0
  \def \uplevel {\begingroup \global \advance \itemdpt by 1 \advance \parindent by 18pt
  \ifcase \itemdpt \or \initem \or \iaitem \or \izitem \fi }
  \def \dnlevel {\par \global \advance \itemdpt by -1 \advance \parindent by -18pt \endgroup }
  \def \itm {\ifcase \itemdpt \or \nitem \or \aitem \or \zitem \fi }
  \def \notext {\vskip -15pt}

  \def \N {{\mathbb {N}}}
  \def \Ninf {{\mathbb {N}\cup\{\infty\}}}
  \def \Z {{\mathbb {Z}}}
  \def \Oo {{\mathcal {O}}}
  \def \OAB {{\mathcal {O}_{A,B}}}
  \def \OGE {{\mathcal {O}_{G,E}}}
  \def \Zplus {\Z ^+}
  \def \OmA {{\Omega_{A}}}
  \def \SZE {{\mathcal {S}_{\Z , E}}}
  \def \GZE {{\cal G}\tight ({\mathcal {S}_{\Z , E}})}


  \def \quoapprox {E^0{\kern -1pt/\kern -2pt}\approx }

  %

  \makeatletter
  \def \Gin@driver{pdftex.def}
  \input color.sty
  \resetatcatcode
  \long \def \Comment #1{\goodbreak \begingroup \color {magenta}\bigskip \hrule \parindent 0pt\parskip 10pt COMMENT: #1
    \par \bigskip \hrule \bigskip \endgroup }

	\def \comment #1{{(\color {magenta} #1})}

	\def \beginRuyComment #1 \endRuyComment {
	  \goodbreak \begingroup \rm \color {magenta}\bigskip \hrule \parindent 0pt\parskip 10pt COMMENT:\par #1
	  \par \bigskip \hrule \bigskip \endgroup \noindent }

	\def \beginEnriqueComment {\goodbreak \begingroup \color {blue}\bigskip \hrule \parindent 0pt\parskip 10pt \par }
	\def \endEnriqueComment {\par \bigskip \hrule \bigskip \endgroup }

  \def\stop{\end}
  \font\fixed = cmitt10 
  \long\def\prelim #1\endprelim{{\overfullrule=0cm \fixed (This is preliminary) #1\par}}
  \def \sep {\medskip \hrule \medskip}
  \def\todo#1{XXXXXXXXX #1 XXXXXXX}

  \def\obsolete{\color{red}}
  \def\endobsolete{\color{black}}

  \def\work{\color{blue}}
  \def\endwork{\color{black}}


\title{A remark on contracting}
\title{inverse semigroups}

\Headlines {A remark on contracting inverse semigroups} {G.~Boava and R.~Exel}

\Authors {Giuliano Boava and Ruy Exel}

\Date {30 August 2014}


\Note {\it Key words and phrases: \rm Inverse semigroup, semi-lattices, tight spectrum, tight filter, ultra-filter,
contracting actions.}

\Note {The first-named author was partially supported by CNPq.}

\bigskip

\Abstract
  A semi-lattice is said to be tree-like when any two of its elements are either orthogonal or comparable.  Given an
inverse semigroup $\S$ whose idempotent semi-lattice is tree-like, and such that all tight filters are ultra-filters, we
present a necessary and sufficient condition for $\S$ to be contracting which looks closer in spirit to the notion of
contracting actions than a recent condition found by the second named author and E.~Pardo.
  \endAbstract

\section Introduction

In a recent paper by the second named author and E.~Pardo \cite{EPFour}, the notion of locally contracting groupoids
introduced by Anantharaman-Delaroche in \cite{AdelaR} was extended to inverse semigroups \cite[Definition 6.4]{EPFour},
as well as to actions of inverse semigroups on topological spaces \cite[Definition 6.2]{EPFour}.  Given an inverse
semigroup $\S$, these concepts were shown to relate to each other via the standard action of $\S$ on $\Eth$, the tight
spectrum of its idempotent semi-lattice.  To be precise it was shown in \cite[Theorem 6.5]{EPFour} that $\S$ is locally
contracting iff\fn
  {The ``if'' part in fact requires that all tight filters be ultra-filters, a condition that has been referred to by the
name of \"{compactable} in \cite{LawsonCompactable}.}
  the standard action of $\S$ on $\Eth$ is locally contracting.

In \cite[Proposition 6.7]{EPFour}, another condition (rephrased here as \ref{MainResult.iii}), which is a lot nicer to
state, and which follows the general paradigm of local contractiveness more closely, was shown to be \"{sufficient} for the
local contractiveness of $\S$.
  In our main result, Theorem \ref{MainResult} below, we take a closer look at this condition and show it to be also
\"{necessary}, provided the inverse semigroup is \"{tree-like}.

One of the key tools to prove our main result is Theorem \ref{CombinResult},  a curious combinatorial fact which  we
suspect may  be  known to specialist in Combinatorial Analysis, but which we have not found anywhere in the literature.

Tree-like inverse semigroups are quite common, especially in the theory of graph C*-algebras (see
e.g.~\cite[Lemma 35.8]{book}), so we feel that this class of inverse semigroups deserves further study.

This  paper is written under  the assumption that the reader is  acquainted with \cite{EPFour}, and to a certain extent
also with \cite{actions}, from where the basic theory of tight representations of inverse semigroups is drawn.

\def\oneton{\{1,\ldots,n\}}

\section A Combinatorial Lemma

In this short section we will prove a crucial combinatorial result to be used in our main result below.

\state Theorem \label CombinResult
  Let $X$ be a set which is decomposed as a finite disjoint union
  $$
  X=\med{\bigsqcup}_{i=1}^n X_i,
  $$
  where each $X_i$ is a nonempty subset,
  and let $f:X\to X$ be a one-to-one map such that
  for every integer $m\geq0,$
  and every $i,j=1,\ldots,n,$   one has that either
  $$  \def\or{,& \hbox{ or}\stake{10pt} \cr }\def\quad{\ }
  \left\{\matrix{
  f^m(X_i)\cap X_j & =         & \emptyset \or
  f^m(X_i)      & \subseteq         & X_j \or
  f^m(X_i)      & \supseteq & X_j.
  } \right.
  \equationmark InclAlternatives
  $$
  Then:
  \izitem
  \zitem
  there exists $i\in\oneton$ and  $m>0,$ such that $f^m(X_i)\subseteq X_i$,
  \zitem
  if  $f$ is not surjective,
  there exists $i\in\oneton$ and  $m>0,$ such that $f^m(X_i)\subsetneqq X_i$.

\Proof Let us begin by proving (i).
  For each $m>0$, let $A^m = \{A^m_{i,j}\}_{i,j}^{\null}$ be the $n\times n$ matrix  defined by
  $$
  A^m_{i,j} =\left\{\matrix{
  1, & \hbox{if } f^m(X_i)\cap X_j \neq\emptyset, \cr \pilar{14pt}
  0, & \hbox{otherwise.}\hfill
  }\right.
  $$

\Case 1: Suppose that there exists some $m>0$, such that every row of $A^m$ has at most one (hence exactly one)
nonzero entry.  In this case, for every $i\in\oneton$, we have that $f^m(X_i)$ intersects a single $X_j$, so it must be contained
in that $X_j$.  We may therefore define a function
  $$
  k:\oneton\to\oneton,
  $$
  such that
  $$
  f^m(X_i) \subseteq X_{k(i)} \for i\in\oneton.
  $$
  From this it easily follows that
  $$
  f^{pm}(X_i) \subseteq X_{k^p(i)}  \for i\in\oneton \for p>0.
  \equationmark PowerContain
  $$

  Given that the set $\{k^p(1):p>0\}$ is finite,
  we may  choose integers $p$ and $q$, with $0<p<q$, and $k^p(1)=k^q(1)$.  Defining
  \def\selectedi{i}
  $$
  \selectedi:=k^p(1)=k^q(1),
  $$
  we then have that  both $f^{pm}(X_1)$ and $f^{qm}(X_1)$ are contained in $X_\selectedi$.
Setting $r = q-p$, observe that 
  $$
  X_\selectedi\supseteq
  f^{qm}(X_1) =
  f^{(r+p)m}(X_1) =
  f^{rm}\big(  f^{pm}(X_1)\big) \subseteq
  f^{rm}(X_\selectedi).
  $$
  The nonempty set   $f^{qm}(X_1)$ is therefore contained in both $X_\selectedi$ and   $f^{rm}(X_\selectedi)$, whence
  $$
  \emptyset\neq  X_\selectedi \cap   f^{rm}(X_\selectedi) \explain\subseteq{PowerContain} X_\selectedi \cap   X_{k^r(i)}.
  $$
  As the $X_j$ are pairwise disjoint, we deduce that $k^r(i)=i$, which is to say that
  $$
  f^{rm}(X_\selectedi) \subseteq X_\selectedi,
  $$
  concluding the proof of (i) in the present case.

\Case 2: Failing the condition characterizing case (1) above,  we are  left with the assumption that, for every $m>0$, at least one  row of $A^m$ has
two or more  nonzero entries.

  Since $A^m$ is a 0-1 matrix,  of which there are only finitely many  (there are exactly $2^{n^2}$ such
matrices), there must be repetitions among the $A^m$, meaning that there are integers  $m_1$ and $m_2$, with $0<m_1<m_2$, and
$A^{m_1}=A^{m_2}$.

Let $p$ be the index of any row of $A^{m_1}$ possessing   two or more nonzero entries, and let $Q$ be the set formed by
the indices
of the columns where such nonzero entries appear, so that $|Q|\geq2$, and
  $$
  A^{m_1}_{p,j} = 1 \iff j\in Q.
  $$
  In particular $f^{m_1}(X_p)$ has a nonempty intersection with $X_q$, for each $q$ in $Q$.   Notice that for any such $q$,
it is impossible that
  $$
  f^{m_1}(X_p) \subseteq X_q,
  $$
  since the $X_i$'s are pairwise disjoint and $f^{m_1}(X_p)$ must intersect at least another $X_i$, given that $|Q|\geq2$.
Thus, given $q$ in $Q$, when comparing $f^{m_1}(X_p)$ with $X_q$ from the point of view of \ref{InclAlternatives}, the
only remaining alternative is  that
  $$
  f^{m_1}(X_p) \supseteq X_q.
  $$

  Since $f^{m_1}(X_p)$ does not intersect $X_j$, for $j\notin Q$, we deduce that
  $$
  f^{m_1}(X_p) = \med\bigcup_{q\in Q} X_q.
  $$

  Observe that $Q$ cannot be equal to $\oneton$, or else
$f^{m_1}(X_p) = X$, and there would be no room for the image of the other $X_i$ under the injective map $f^{m_1}$.
Consequently
  $$
  2\leq|Q|<n.
  $$

  Recalling that $A^{m_1}=A^{m_2}$, the above argument also proves that
  $$
  f^{m_2}(X_p) = \med\bigcup_{q\in Q} X_q.
  $$
  Defining   $m:=m_2-m_1$, notice that
  $$
  f^m\Big(\med\bigcup_{q\in Q} X_q\Big) =
  f^m\big(f^{m_1}(X_p)\big) =
  f^{m+m_1}(X_p) =    f^{m_2}(X_p) =
  \med\bigcup_{q\in Q} X_q,
  $$
  so we may restrict  $f^m$ to
  $$
  X' := \med\bigcup_{q\in Q} X_q
  $$
  obtaining an injective map
  $$
  f^m:X'\to X',
  $$
  satisfying all of the assumptions of the statement, with $X'$ decomposing into a smaller number of parwise disjoint
components $X_q$'s.  Therefore the result follows immediately by induction on $n$.

In order to prove (ii) we may now use (i) and hence we may assume that
  $$
  f^m(X_i)\subseteq X_i,
  $$
  for some $i\in\oneton$, and some $m>0$.  In case the above is a proper inclusion we are done, so we assume the contrary,
meaning that $f^m(X_i)=X_i$.  Setting
  $$
  X'= \med\bigcup_{j\neq i}X_j,
  $$
  and recalling that $f^m$ is injective, we then have that
  $$
  f^m(X')\subseteq X'.
  $$
  Observe that this must is a proper
inclusion since  otherwise $f^m$ would be surjective, contradicting the hypothesis. Therefore  the conclusion  follows
again by
induction on $n$.
  \endProof

\def\<{\ll}
\def\Dmt#1{\Dm \theta {#1}}

\section Preliminaries on semi-lattices and inverse semigroups

In this section we will freely use the notation introduced in \cite{EPFour}.  Our main goal will be to
introduce the class of inverse semigroups to which our main result applies.  We will  also prove
some basic related results.

\definition
  A semi-lattice $\E$ with zero is called \"{tree-like} if, for any $e$ and $f$ in $\E$, one has that
  $$
  e\perp f,\quad e\leq f,\quad \hbox{or}\quad e\geq f.
  $$
  If $\S$ is an inverse semigroup whose idempotent semi-lattice is tree-like, we will say that $\S$ is \"{tree-like}.

Given  $e$ and $f$ in a semi-lattice, it is easy to see that
  $$
  e\leq f \imply \Dmt e \subseteq \Dmt f.
  \equationmark InclusDomains
  $$
  The converse of this fact is however not  true.  For example, if  $\S$ is obtained by adding  a zero
element to an inverse semigroup without a zero, then
  $$
  \xi:=\E\setminus \{0\}
  $$
  is the only ultra-filter on $\E$.  Consequently
$\Eth=\{\xi\}$, and then $\Dmt e = \{\xi\}$ for every nonzero idempotent $e$, so the converse of \ref{InclusDomains} is seen to
fail badly.
However, in a tree-like inverse semigroup, it is  easy to see that
  $$
  \Dmt e \subsetneqq \Dmt f \imply e\leq f,
  \equationmark TreeLikeDomains
  $$
  simply because,  under the  assumption that $\Dmt e \subsetneqq \Dmt f$, alternatives ``$e\perp f$'' and  ``$e\geq f$'' are clearly excluded.

The strict inclusion above has other interesting consequences.  Assuming that    $\Dmt f \setminus \Dmt e$ is indeed
nonempty, and noticing that it  is an open subset of $\Eth$, we may find an  ultra-filter $\xi$ there, meaning that  $f\in\xi$ and $e\notin\xi$.
By \cite[Lemma 12.3]{actions}, it follows that there exists some $d\in\xi$ such that $d\perp e$, and upon replacing $d$ with $df$,
we may clearly assume that $d\leq f$.

\definition Given $e$ and $f$ in a semi-lattice $\E$, we will say that $e\<f$, whenever  $e\leq f$, and there exists a nonzero $d\leq f$
such that $d\perp e$.

Using this terminology we may then state the following  fact:

\state Proposition  \label StrongInclusion
  Let $\E$ be a tree-like semi-lattice.  Then
  $$
  \Dmt e \subsetneqq \Dmt f \imply e\<f \for e,f\in\E.
  $$

There is another slightly annoying question related to the converse of \ref{InclusDomains} which we would like to get
out of our way as soon as possible:

\state Proposition \label Anoying
  Given any inverse semigroup $\S$, and given $s\in\S$, and $e\in\E$, we have that:
  \izitem
  \zitem if $e\leq s^*s$, then $\Dmt e$ is contained in the domain of\/ $\theta_s$, and $\theta_s(\Dmt e) = \Dmt {ses^*}$.
  \zitem if $\Dmt e$ is contained in the domain of\/ $\theta_s$, then $\theta_s(\Dmt e) = \Dmt {ses^*}$.

\Proof
  The first assertion in (i) is obvious.  As for the second, recall that $\theta_e$ is the identity map on $\Dmt e$ so, in
particular, the range of $\theta_e$ is $\Dmt e$.  Thus,
  $$
  \theta_s(\Dmt e) =
  \theta_s\big(\hbox{Ran}(\theta_e)\big) =
  \hbox{Ran}(\theta_s\theta_e) =
  \hbox{Ran}(\theta_{se}) =
  \Dmt  {se (se)^*} =
  \Dmt  {ses^*},
  $$
  proving (i).

As the reader may have already anticipated, the catch in (ii) is that it is  assumed that $\Dmt e\subseteq\Dmt {s^*s}$, but not
necessarily that $e\leq s^*s$.  Fortunately this can be easily  circumvented as follows:
  $$
  \theta_s(\Dmt e) =
  \theta_s(\Dmt e\cap\Dmt {s^*s}) =
  \theta_s(\Dmt {es^*s}) \explica={(i)}
  \Dmt  {s(es^*s)s^*} =
  \Dmt  {ses^*}.
  \endProof

\section The main result

Given the above preparations, we are now ready to prove our main result.

\state Theorem \label MainResult
  Let $\S$ be a tree-like  inverse semigroup such that every tight filter in $\E$ is an ultra-filter.  Then the
following are equivalent:
  \izitem
  \zitem $\S$ is locally contracting,
  \zitem the standard action $\theta :\S \curvearrowright \Eth$ is locally  contracting,
  \zitem for every nonzero $e\in\E$, there exists an idempotent  $f\leq e$, and an element $s\in\S$, such that $f\leq s^*s$,
and $sfs^*\<f$.

\Proof  The equivalence between (i) and (ii) follows from \cite[Theorem 6.5]{EPFour}.

In order to prove that (iii) implies (i), given a nonzero $e$ in $\E$, let $f$ and $s$ be as in (iii).  Since $sfs^*\<f$,
there exists a nonzero $f_0\leq f$, such that $f_0\perp sfs^*$, and then we see that $f_0$ together with $f_1:=f$ obey the
conditions of \cite[Proposition 6.7]{EPFour}, from where one  deduces that $\S$ is locally contracting, proving (i).

The most delicate part of this proof is the implication (i)$\Rightarrow$(iii), which we take up next.
  Given a nonzero $e\in\E$, let $U=\Dmt e$, and choose $s$ and $V$ as in \cite[Definition 6.2]{EPFour}
  so that,
  $$
  V\subseteq U, \quad \clos V\subseteq\Dmt  {s^*s} \and  \theta_s(\clos V) \subsetneqq V.
  $$

In the first part of the proof we will show that $V$ may be chosen to be of the form
  $$
  V=\medcup_{f\in F} \Dmt  {f},
  $$
  where $F$ is a finite set of idempotents satisfying $f\leq es^*s$.

  In order to achieve this,   for each $\xi$ in $\theta_s(\clos V)$, choose a neighborhood of $\xi$ contained in $V$.  By hypothesis we have that $\xi$ is an
ultra-filter, and by
  \cite[Proposition 2.5]{EPFour} we may suppose that such a neighborhood is of the form $\Dmt {f_\xi}$, for some $f_\xi$ in $\E$,
whence
  $$
  \xi\in\Dmt  {f_\xi}\subseteq V,
  $$
  so we see that   the $\Dmt  {f_\xi}$  form an open cover for $\theta_s(\clos V)$.
  Since
  $$
  V\subseteq U\cap\Dmt  {s^*s} =\Dmt  e\cap\Dmt  {s^*s} = \Dmt  {es^*s},
  \equationmark VInEssStar
  $$
  then also each $\Dmt {f_\xi}\subseteq\Dmt {es^*s}$, and therefore
  $$
  \Dmt  {f_\xi} =
  \Dmt  {es^*s} \cap \Dmt  {f_\xi} =
  \Dmt  {es^*sf_\xi}.
  $$

  Upon  replacing  each $f_\xi$  by
  $$
  f_\xi':= es^*sf_\xi,
  $$
  we may therefore assume that $f_\xi\leq es^*s$.

  Being a closed subset of $\Dmt  {s^*s}$, observe that $\clos V$ is   compact, and hence so  is $\theta_s(\clos
V)$.  We may then take a finite subcover of the above cover, say
  $$\thickmuskip 12mu
  \theta_s(\clos V)\subseteq\medcup_{f\in F'} \Dmt  {f},
  \equationmark LocalContrEquivClosV
  $$
  where $F'$ is a finite set consisting of some of the $f_\xi$.

We next claim that there exists  a nonzero idempotent $f_0\leq es^*s$, such that
  $$
  \Dmt  {f_0}\subseteq V\setminus \theta_s(\clos V).
  \equationmark LocalContrEquivFZero
  $$

  To see this, first observe that
  $
  V\setminus \theta_s(\clos V)
  $
  is  open and nonempty.
Even  without assuming that all tight filters are ultra-filters, we may use the density of the set formed by the latter
to find some ultra-filter $\xi$ in   $V\setminus \theta_s(\clos V)$.  An application of
  \cite[Proposition 2.5]{EPFour}
  then provides $f_0$ in $\E$ such that
  $$
  \xi \subseteq \Dmt  {f_0} \subseteq V\setminus \theta_s(\clos V),
  \equationmark FZeroAppears
  $$
  and, again by  \ref{VInEssStar}, we may assume that $f_0\leq es^*s$.  Adding $f_0$  to  $F'$, we form the set
  $$
  F := \{f_0\}\cup F',
  $$
  with which we define
  $$
  W:=\medcup_{f\in F} \Dmt  {f}.
  \equationmark definitionOfW
  $$

  We then have that  $W$ is clopen, and that
  $$
  \theta_s(\clos V)\subsetneqq W\subseteq V,
  $$
  where the proper inclusion above is a consequence of \ref{FZeroAppears} and
the fact that we have included $f_0$ in $F$.  Applying $\theta_s$ to the sets above we then deduce that
  $$
  \thickmuskip 12mu
  \theta_s(W) \subseteq \theta_s(V) \subseteq \theta_s(\clos V) \subsetneqq W.
  $$
  This completes the task outlined at the beginning of the proof.

Notice that for any  $f_1$ and $f_2$ in $F$, we have  by (i) that
  $$
  f_1\perp f_2, \quad f_1\leq f_2, \quad\hbox{or}\quad  f_1\geq f_2,
  $$
  in which case
  $$
  \Dmt {f_1}\cap\Dmt {f_2}=\emptyset, \quad \Dmt {f_1}\subseteq\Dmt {f_2}, \quad \hbox{or} \quad  \Dmt {f_1}\supseteq\Dmt {f_2},
  $$
  respectively.  Replacing $F$ by the subset of its maximal elements we may then assume that $F$ is formed by  pairwise
orthogonal idempotents, in which case  \ref{definitionOfW} is a disjoint union.

In order to proceed, let us consider two cases.

\Case 1:  Assuming that  $|F|=1$, say $F=\{f\}$, we have that $W=\Dmt f$, and then
  $$
  \Dmt {sfs^*} \={Anoying.i} \theta_s(\Dmt f) = \theta_s(W) \subsetneqq W = \Dmt f,
  $$
  so   \ref{StrongInclusion} implies that $sfs^*\< f$, concluding the proof.

\Case 2:  Assuming that  $|F|>1$, let us  consider $\theta_s$ as a function
  $$
  \theta_s:W\to W,
  $$
  observing that it is an injective but not surjective map.  Using \ref{CombinResult} and the fact that \ref{definitionOfW} is a
disjoint union, we have that
  $$
  \theta_s^m(\Dmt f)\subsetneqq\Dmt f,
  $$
  for some integer $m>0$, and some $f$ in $F$.  We next notice that
  $$
  \theta_s^m(\Dmt f) =   \theta_{s^m}(\Dmt f) \={Anoying.ii}  \Dmt {s^mfs^{*m}},
  $$
  so we deduce from the above that
  $$
  \Dmt {s^mfs^{*m}} \subsetneqq\Dmt f,
  $$
  and then  \ref{StrongInclusion} implies that $s^mfs^{*m}\< f$.

To conclude we must still address the requirement that $f\leq s^{*m}s^m $. For  this  observe that since $W$
is contained in the domain of $\theta_s$, and since $W$ is invariant under $\theta_s$, we have that $W$ is also contained in the domain
of  $\theta_s^m$, namely
  $$
  W\subseteq\Dmt {s^{*m} s^m}.
  $$

  Recalling that we are working under the hypothesis that $|F|>1$, we have that
  $\Dmt f$ is a proper subset of $W$, so
  $$
  \Dmt f \subsetneqq W \subseteq \Dmt {s^{*m} s^m},
  $$
  and then
  $f \< s^{*m} s^m$, by \ref{StrongInclusion}.
\endProof

The main point we would like to make in the present work is that, even though the definition of locally contracting
actions given in
  \cite[Definition 6.2]{EPFour}
  is syntactically closer to
  \ref{MainResult.iii},
  from a logical point of view, the former  is closer to the
  notion of locally contracting inverse semigroup described in
  \cite[Definition 6.4]{EPFour},
  since these are  equivalent to each other under broader conditions, as proved in
  \cite[Theorem 6.5]{EPFour}

As seen in Theorem \ref{MainResult}, above, all of these are equivalent to each other under the rather strong assumption that $\S$
is tree-like, but it would be highly desirable to decide if the tree-like property is indeed necessary for the proof of
\ref{MainResult}.

A related  problem  is to decide conditions under which  the converse of \cite[Proposition 6.3]{EPFour} also holds.

\references

\Article AdelaR
  C. Anantharaman-Delaroche;
  Purely infinite $C^*$-algebras arising form dynamical systems;
  Bull. Soc. Math. France, 125 (1997), no. 2, 199-225

\Article actions
  R. Exel;
  Inverse semigroups and combinatorial C*-algebras;
  Bull. Braz. Math. Soc., 39 (2008), no. 2, 191-313

\Bibitem book
  R. Exel;
  Partial Dynamical Systems, Fell Bundles and Applications;
  {Licensed under a Creative Commons Attribution-ShareAlike 4.0 International License, 351pp, 2014.  Available online from\break
mtm.ufsc.br/$\sim$exel/publications. \rm PDF file md5sum: bc4cbce3debdb584ca226176b9b76924}

\Bibitem EPFour
  R. Exel and E. Pardo;
  The tight groupoid of an inverse semigroup;
  arXiv:1408.5278 [math.OA], 2014

\Article LawsonCompactable
  M. V. Lawson;
  Compactable semilattices;
  Semigroup Forum, 81 (2010), no. 1, 187-199

\endgroup

\bigskip \noindent \tensc
  Departamento de Matem\'atica;
  Universidade Federal de Santa Catarina;
  88010-970 Florian\'opolis SC;
  Brazil


  \close

\bye